# A Condition for Hopf bifurcation to occur in Equations of Lotka-Volterra Type with Delays

Chol Kim

Faculty of Mathematics, **Kim Il Sung** University, Pyongyang, DPRK

**Abstract**

It is known that Lotka-Volterra type differential equations with delays or distributed delays have an important role in modeling ecological systems. In this paper we study the effects of distributed delay on the dynamics of the harvested one predator − two prey model. Using the expectation of the distribution of the delay as a bifurcation parameter, we show that the equilibrium that was asymptotic stable becomes unstable and Hopf bifurcation can occur as the expectation crosses some critical values.

**Keyword:** Lotka-Volterra equation, distributed delay, Hopf bifurcation, stability, equilibrium

**2010 MSC:** 34K18, 37G15

## 1. Introduction

In this paper we study the stability and bifurcation of the following one predator − two prey model with distributed delay.

$$\begin{cases} x_1'(t) = x_1(t)(r_1 - a_{11}x_1(t) - a_{12}x_2(t) - a_{13}x_3(t)) \\ x_2'(t) = x_2(t)(r_2 - a_{21}x_1(t) - a_{22}x_2(t) - a_{23}x_3(t)) \\ x_3'(t) = x_3(t)\left(-r_3 + a_{31}\int_0^\infty f_E(\tau)x_1(t-\tau)d\tau + a_{32}\int_0^\infty f_E(\tau)x_2(t-\tau)d\tau\right) - H \end{cases} \quad (1)$$

The meaning of $x_i(t)$, $a_{ij}$, $r_i$, $H$ ($i, j = 1, 2, 3$) is the same as those in paper [1].

Kumar et al. [5] considered the following one predator − two prey system without delay.

$$\begin{cases} x_1'(t) = x_1(t)(r_1 - a_{11}x_1(t) - a_{12}x_2(t) - a_{13}x_3(t)) \\ x_2'(t) = x_2(t)(r_2 - a_{21}x_1(t) - a_{22}x_2(t) - a_{23}x_3(t)) \\ x_3'(t) = x_3(t)(-r_3 + a_{31}x_1(t) + a_{32}x_2(t)) - H \end{cases} \quad (2)$$

They studied the stability and bifurcation of system (2) using the harvesting rate $H$ as a bifurcation parameter.

Zhihua Liu et al. [1] studied the stability of positive equilibrium and Hopf bifurcation of system (1) using the delay value $\tau$ as a bifurcation parameter when $f_E(x) = \delta(x-\tau)$ ($\delta$-function).

Systems with delays or distributed delays are studied in many situations (see [2, 3, 4, 6, 7, 8, 9]). In [7] the linear differential equation with distributed delay

$$x'(t) = -\alpha x(t) - \beta \int_0^\infty f(\tau) x(t-\tau) d\tau \quad (3)$$

is concerned. They studied sufficient conditions for stability of equilibrium of system (3).

In this paper, we study the effect of distributed delays on system (1).

## 2. Main result

Like Kumar et al. [5], we also choose $r_1 = r_2 = r_3 = a_{11} = a_{22} = a_{12} = a_{23} = 1$, $a_{33} = 0$, $a_{21} = 1.5$, $a_{32} = 0.5$, $a_{13} = a$, $a_{31} = \frac{1}{2}a$. Then system (1) becomes as following.

$$\begin{cases} x_1'(t) = x_1(t)(1 - x_1(t) - x_2(t) - ax_3(t)) \\ x_2'(t) = x_2(t)(1 - 1.5x_1(t) - x_2(t) - x_3(t)) \\ x_3'(t) = x_3(t)\left(-1 + a/2\int_0^\infty f_E(\tau) x_1(t-\tau) d\tau + 1/2\int_0^\infty f_E(\tau) x_2(t-\tau) d\tau\right) - H \end{cases} \quad (4)$$

Here $f_E(\tau)$ satisfies the following conditions.

1) $f_E(\tau)\Delta\tau$ represents the probability of the event of occurrence of delay between $\tau$ and $\tau+\Delta\tau$. That is, $f_E(\tau)$ satisfies the following:

$$\forall \tau \in [0,\infty): f_E(\tau) \geq 0, \quad \int_0^\infty f_E(\tau)d\tau = 1, \quad \int_0^\infty \tau f_E(\tau)d\tau = E.$$

Here we call $E$ the *expectation* of *distributed delay*.

2) As the variable $\tau$ tends infinity, it is rapidly decreasing or it's support is compact.

Let's denote the positive equilibrium of (4) with $E^* = (x_1^*, x_2^*, x_3^*)$. The $x_i^*$ is the same as paper [5], considering the condition (1) of $f_E(\tau)$.

$$x_1^* = 2(a-1)x_3^*,$$

$$x_2^* = 1 - (3a-2)x_3^*, \quad (5)$$

$$x_3^* = \frac{\left(1+\sqrt{1+8H(a-2)(2a-1)}\right)}{2(a-2)(2a-1)}$$

Here $H > 0$ and $x_i^* > 0 \ (i = \overline{1,3})$, thus we assume that

$$H < H_c = \frac{a^2 - 4a + 2}{3a - 2}, \quad a > 2 + \sqrt{2}.$$

If there is no distributed delay, the system (4) is the same with the system

$$\begin{cases} x_1'(t) = x_1(t)(1 - x_1(t) - x_2(t) - ax_3(t)) \\ x_2'(t) = x_2(t)(1 - 1.5x_1(t) - x_2(t) - x_3(t)) \\ x_3'(t) = x_3(t)(-1 + a/2\, x_1(t) + 1/2\, x_2(t)) - H \end{cases} \quad (6)$$

which has been considered in [5]. The positive equilibrium of (6) is equal to the positive equilibrium (5) of (4).

Hence if there is no effect of distributed delay, that is, when $E = 0$, the positive equilibrium of (4) is locally asymptotically stable if and only if the following conditions hold [1]:

$$a_1 > 0, \quad a_2 + a_4 < 0, \quad a_3 + a_5 > 0, \quad a_1(a_2 + a_4) > a_3 + a_5. \tag{7}$$

Here

$$a_1 = x_1^* + x_2^* - \frac{H}{x_3^*}, \qquad a_2 = -\frac{1}{2}x_1^* x_2^* - \left(x_1^* + x_2^*\right)\frac{H}{x_3^*},$$

$$a_3 = \frac{1}{4}x_1^* x_2^* x_3^* (2a-1)(a-2), \quad a_4 = \frac{1}{2}\left(a^2 x_1^* + x_2^*\right)x_3^*, \quad a_5 = \frac{1}{2}x_1^* x_2^* \frac{H}{x_3^*}.$$

Now we study the effect of distributed delays on system (4).

Let $N_i(t) = x_i(t) - x_i^*$, $(i = 1, 2, 3)$, the system (4) becomes

$$N_1'(t) = -x_1^*\big(N_1(t) + N_2(t) + aN_3(t)\big) - N_1(t)\big(N_1(t) + N_2(t) + aN_3(t)\big),$$

$$N_2'(t) = -x_2^*\left(\frac{3}{2}N_1(t) + N_2(t) + N_3(t)\right) - N_2(t)\left(\frac{3}{2}N_1(t) + N_2(t) + N_3(t)\right),$$

$$N_3'(t) = x_3^*\left(\frac{a}{2}\int_0^\infty f_E(\tau) N_1(t-\tau)d\tau + \frac{1}{2}\int_0^\infty f_E(\tau) N_2(t-\tau)d\tau + \frac{H}{\left(x_3^*\right)^2}N_3(t)\right) +$$

$$+ N_3(t)\left(\frac{a}{2}\int_0^\infty f_E(\tau) N_1(t-\tau)d\tau + \frac{1}{2}\int_0^\infty f_E(\tau) N_2(t-\tau)d\tau\right).$$

The characteristic equation of the linearized system is as follows:

$$\lambda^3 + a_1 \lambda^2 + a_2 \lambda + a_3 \int_0^\infty f_E(\tau) e^{-\lambda \tau} d\tau + a_4 \lambda \int_0^\infty f_E(\tau) e^{-\lambda \tau} d\tau + a_5 = 0. \tag{9}$$

**Theorem 1.** *Suppose the conditions in* (7) *hold. Then the equation* (9) *has a simple pair of conjugate purely imaginary roots* $\pm i\omega_1$ *with some expectation* $E_1$.

**Proof:** $\lambda = 0$ is a root of (9) if and only if $a_3 + a_5 = 0$. By the condition of the theorem, we have $a_3 + a_5 > 0$ and hence $\lambda = 0$ is not a root of (9).

Assume that for some expectation $E$, $\lambda = i\omega_1$ $(\omega_1 > 0)$ is a root of (9). Substituting $\lambda = i\omega_1$ into (9) and separating the real and imaginary parts yield

$$a_3 \int_0^\infty f_E(\tau)\cos(\omega_1 \tau)d\tau + a_4 \omega_1 \int_0^\infty f_E(\tau)\sin(\omega_1 \tau)d\tau = a_1 \omega_1^2 - a_5, \tag{10}$$

$$-a_3\int_0^\infty f_E(\tau)\sin(\omega_1\tau)d\tau + a_4\omega_1\int_0^\infty f_E(\tau)\cos(\omega_1\tau)d\tau = \omega_1^3 - a_2\omega_1. \qquad (11)$$

Adding sidewise after squaring the left and right sides of (10) and (11), we can obtain $F(\omega_1) = G(\omega_1)$. Here

$$F(\omega) := \left(a_3\int_0^\infty f_E(\tau)\cos(\omega\tau)d\tau + a_4\omega\int_0^\infty f_E(\tau)\sin(\omega\tau)d\tau\right)^2 +$$

$$+ \left(-a_3\int_0^\infty f_E(\tau)\sin(\omega\tau)d\tau + a_4\omega\int_0^\infty f_E(\tau)\cos(\omega\tau)d\tau\right)^2 =$$

$$= a_3^2\left[\left(\int_0^\infty f_E(\tau)\cos(\omega\tau)d\tau\right)^2 + \left(\int_0^\infty f_E(\tau)\sin(\omega\tau)d\tau\right)^2\right] +$$

$$+ a_4^2\omega^2\left[\left(\int_0^\infty f_E(\tau)\sin(\omega\tau)d\tau\right)^2 + \left(\int_0^\infty f_E(\tau)\cos(\omega\tau)d\tau\right)^2\right],$$

$$G(\omega) := (a_1\omega^2 - a_5)^2 + (\omega^3 - a_2\omega)^2.$$

Here using the Schwartz inequalities, we have

$$\left(\int_0^\infty f(\tau)\cos(\omega\tau)d\tau\right)^2 \leq \int_0^\infty f(\tau)d\tau \int_0^\infty f(\tau)\cos^2(\omega\tau)d\tau,$$

$$\left(\int_0^\infty f(\tau)\sin(\omega\tau)d\tau\right)^2 \leq \int_0^\infty f(\tau)d\tau \int_0^\infty f(\tau)\sin^2(\omega\tau)d\tau.$$

Thus we have

$$\forall \omega; \; F(\omega) \leq a_3^2 + a_4^2\omega^2.$$

Let's consider $H(\omega) := F(\omega) - G(\omega)$. Since $F, G \in C(R)$, hence $H \in C(R)$. If we let $\omega = 0$, then $H(0) = F(0) - G(0) = a_3^2 - a_5^2 > 0$. (The fact that $a_3^2 - a_5^2 > 0$ has been proved in the proof of the theorem 1 of [1].) Now if we let $\omega = \omega_0$ where $\omega_0$ is equal to $\sigma_+$ in [1], then we have

$$H(\omega_0) = F(\omega_0) - G(\omega_0) = F(\omega_0) - (a_3^2 + a_4^2 \omega_0^2) \leq 0.$$

(The fact that $G(\omega_0) = (a_3^2 + a_4^2 \omega_0^2)$ has been proved in the proof of the theorem 1 of [1].) Thus there exists a $\omega_1 \in (0, \omega_0]$ such that $F(\omega_1) = G(\omega_1)$.

From (10) and (11),

$$\int_0^\infty f_E(\tau)\cos(\omega_1 \tau)d\tau = \frac{\omega_1(a_4\omega_1^2 - a_2 a_4 + a_1 a_3) - a_3 a_5}{a_3^2 + a_4^2 \omega_1^2}, \qquad (12)$$

On the other hand, in [7] we know

$$\int_0^\infty f_{E_1}(\tau)\cos(\omega_1 \tau)d\tau \geq 1 - \frac{c\omega_1}{\pi}E_1, \text{ where } c := \sup\left\{c \middle| \cos x = 1 - \frac{cx}{\pi}, x > 0\right\} \approx 2.2764.$$

By (12),

$$E_1 \geq \frac{\pi[a_3^2 + a_4^2 \omega_1^2 + a_3 a_5 - \omega_1^2(a_4 \omega_1^2 - a_2 a_4 + a_1 a_3)]}{c\omega_1(a_3^2 + a_4^2 \omega_1^2)}.$$

The equation (9) has roots $\pm i\omega_1$ when the expectation takes the value $E_1$ such that the above inequality holds.(QED)

**Theorem 2.** *Suppose the conditions in* (7) *hold. Then we have*

1) *If* $E \in [0, E_1)$, *then the positive equilibrium* $E^*$ *of* (4) *is locally asymptotically stable*.

2) *If the condition*

$$\text{Re}\left\{-\left.\frac{3\lambda^2+2a_1\lambda+a_2+a_4\int_0^\infty f_E(\tau)e^{-\lambda\tau}d\tau-a_3\int_0^\infty \tau f_E(\tau)e^{-\lambda\tau}d\tau-a_4\lambda\int_0^\infty \tau f_E(\tau)e^{-\lambda\tau}d\tau}{a_3\int_0^\infty \frac{df_E(\tau)}{dE}e^{-\lambda\tau}d\tau+a_4\lambda\int_0^\infty \frac{df_E(\tau)}{dE}e^{-\lambda\tau}d\tau}\right|_{\lambda=i\omega_1}\right\}\neq 0$$

holds, then in the system (4) *Hopf bifurcation occurs when* $E=E_1$.

**Proof:** If $E=0$, the characteristic equation (9) is as follows:

$$\lambda^3+a_1\lambda^2+(a_2+a_4)\lambda+a_3+a_5=0.$$

(See the theorem 4.0.5 of [7].) This equation has roots with negative real part. If $E=E_1$, the equation (9) has a pair of conjugate purely imaginary roots $\pm i\omega_1$. If $E\in[0,E_1)$, the equation (9) has roots with negative real part. Thus if $E\in[0,E_1)$, the positive equilibrium $E^*$ of (4) is locally asymptotically stable (4).

Let denote the root of (9) by $\lambda(E)=\mu(E)+i\omega(E)$, then we have

$$\mu(E_1)=0, \omega(E_1)=\omega_1.$$

If the transversal condition $\mu'(E_1)=\left.\dfrac{d\operatorname{Re}\lambda(E)}{dE}\right|_{E=E_1}\neq 0$ is satisfied, then the Hopf bifurcation occurs in the system (4) when $E=E_1$. If we denote (9) by $G(\lambda,E)=0$, then

$$\left.\frac{dE}{d\lambda}\right|_{\lambda=i\omega_1}=-\left.\frac{G_\lambda}{G_E}\right|_{\lambda=i\omega_1}=$$

$$=-\left.\frac{3\lambda^2+2a_1\lambda+a_2+a_4\int_0^\infty f_E(\tau)e^{-\lambda\tau}d\tau-a_3\int_0^\infty \tau f_E(\tau)e^{-\lambda\tau}d\tau-a_4\lambda\int_0^\infty \tau f_E(\tau)e^{-\lambda\tau}d\tau}{a_3\int_0^\infty \frac{df_E(\tau)}{dE}e^{-\lambda\tau}d\tau+a_4\lambda\int_0^\infty \frac{df_E(\tau)}{dE}e^{-\lambda\tau}d\tau}\right|_{\lambda=i\omega_1}.$$

Thus if the condition 2) of the theorem holds, the Hopf bifurcation occurs as the transversal condition holds. (QED)

## 3. Conclusion and Further Study

In this paper we provided a condition for Hopf bifurcation to occur in the equilibrium in equations of Lotka-Volterra type with distributed delays when using the expectation of the distribution of the delay as a bifurcation parameter. (See (7) and theorem 1, 2.)

Lotka-Volterra competitive systems (with delay) are being studied in recent, too. For example, [10] proposed a discrete Lotka-Volterra competition system with infinite delays and feedback controls. Sufficient conditions which ensure the global attractivity of the system are obtained. In [12], a discrete nonautonomous two-species Lotka-Volterra competitive system with delays and feedback controls is proposed and investigated. By using the method of discrete Lyapunov functionals, new sufficient conditions on the permanence of species and global attractivity of the system are established. [11] discussed a two-species discrete competition system. The local stability of positive equilibrium is obtained. Further, a sufficient condition for the global asymptotic stability of positive equilibrium is established.

In further study we will study bifurcation problems and control problems for Lotka-Volterra type systems with various type of delays including infinite delay.

## 4. ACKNOWLEDGMENTS

Authors would like to thank anonymous reviewers' help and advice.